\documentclass[12pt]{article}

\usepackage{amssymb,amsmath,exscale}

\usepackage{graphicx}

\usepackage{color}
\definecolor{marin}{rgb}   {0.,   0.3,   0.7} 
\definecolor{rouge}{rgb}   {0.8,   0.,   0.} 
\definecolor{sepia}{rgb}   {0.8,   0.5,   0.} 
\usepackage[colorlinks,citecolor=marin,linkcolor=rouge,
            bookmarksopen,
            bookmarksnumbered
           ]{hyperref}

\newtheorem{lemma}{Lemma}[section]

\newtheorem{theorem}[lemma]{Theorem}
\newtheorem{proposition}[lemma]{Proposition}

\newtheorem{remark}[lemma]{Remark}
\newtheorem{example}[lemma]{Example}
\newtheorem{hypothesis}[lemma]{Hypothesis}
\newtheorem{notation}[lemma]{Notation}
\newtheorem{definition}[lemma]{Definition}
\newtheorem{conclusion}[lemma]{Conclusion}
\numberwithin{equation}{section}
\newcommand{\QED}{\mbox{}\hfill \raisebox{-0.2pt}{\rule{5.6pt}{6pt}\rule{0pt}{0pt}} 
          \medskip\par}             
\newenvironment{Proof}{\noindent
    \parindent=0pt\abovedisplayskip = 0.5\abovedisplayskip 
    \belowdisplayskip=\abovedisplayskip{\bfseries Proof. }}{\QED}
\newenvironment{Proofof}[1]{\noindent
    \parindent=0pt\abovedisplayskip = 0.5\abovedisplayskip
    \belowdisplayskip=\abovedisplayskip{\bfseries Proof of #1. }}{\QED}

\newcommand{\Ac}{\mathcal{A}}

\newcommand{\dd}{\mathrm{d}}

\newcommand{\Ic}{\mathcal{I}}

\newcommand{\jb}{{\boldsymbol{j}}}
\newcommand{\kb}{{\boldsymbol{k}}}

\newcommand{\N}{\mathbb{N}}
\newcommand{\Nc}{\mathcal{N}}

\newcommand{\R}{\mathbb{R}}

\newcommand{\C}{\mathbb{C}}

\newcommand{\T}{\mathbb{T}}
\newcommand{\Z}{\mathbb{Z}}
\newcommand{\Tc}{\mathcal{T}}

\newcommand{\Uc}{\mathcal{U}}

\newcommand{\Wc}{\mathcal{W}}

\newcommand{\Zc}{\mathcal{Z}}

\newcommand{\Norm}[2]{\|#1\|\left.\vphantom{T_{j_0}^0}\!\!\right._{#2}}         
\newcommand{\SNorm}[2]{|#1|\left.\vphantom{T_{j_0}^0}\!\!\right._{#2}}             

\title{Birkhoff normal form and splitting methods\\ for semi linear Hamiltonian PDEs. \\
Part I: Finite dimensional discretization.}        
\author{Erwan Faou, Beno\^it Gr\'ebert and Eric Paturel}       
\begin{document}
\maketitle
\abstract{
We consider {\em discretized}  Hamiltonian PDEs associated with a Hamiltonian function  that can be split into a linear unbounded operator and a regular nonlinear part. We consider splitting methods associated with this decomposition. Using a finite dimensional Birkhoff normal form result, we show the almost preservation of the {\em actions} of the numerical solution associated with the splitting method over arbitrary long time, provided the Sobolev norms of the initial data is small enough, and for asymptotically large level of space approximation. This result holds under {\em generic}  non resonance conditions on the frequencies of the linear operator and on the step size. We apply this results to nonlinear Schr\"odinger equations as well as the nonlinear wave equation. }

\tableofcontents

\section{Introduction}

In this work, we consider a class of Hamiltonian partial differential equations whose Hamiltonian functionals $H=H_0+P$ can be divided into a linear unbounded operator $H_0$ with discrete spectrum and a nonlinear function $P$ having a zero of order at least $3$ at the origin of the phase space. Typical examples are given by the nonlinear wave equation  or the nonlinear Schr\"odinger equation on the torus. We consider discretizations of this PDEs and denote by $H^{(K)}=H_0^{(K)}+P^{(K)}$ the corresponding discrete Hamiltonian, where $K$ is a discretization parameter. Typically, $K$ denotes a spectral parameter in a collocation method. 

Amongst all the numerical schemes that can be applied to these Hamiltonian PDEs, splitting methods entail many advantages, as they provide symplectic and explicit schemes, and can be easily implemented using fast Fourier transform if the spectrum of $H_0$ expresses easily in  Fourier basis. Generally speaking, a splitting scheme is based on the approximation 
\begin{equation}
\label{E001}
\varphi_{H^{(K)}}^h \simeq \varphi_{H_0^{(K)}}^h \circ \varphi_{P^{(K)}}^h
\end{equation}
for small time $h$, and where $\varphi_{Q}^t$ denotes the exact flow of the Hamiltonian system associated with the Hamiltonian function $Q$. For a given time $t = nh$, $n \in \N$, the solution starting at some initial value $z^0$ is then approximated by 
\begin{equation}
\label{Ezn}
\varphi_{H^{(K)}}^t(z^0)  \simeq z^n = \Big(\varphi_{H_0^{(K)}}^h \circ \varphi_{P^{(K)}}^h\Big)^n (z^0). 
\end{equation}

The understanding of the long-time behavior of splitting methods for Hamiltonian PDEs is a fundamental ongoing challenge in the field of geometric integration, as the classical arguments of {\em backward error analysis} (see for instance \cite{HLW}) do not apply in this situation, where the frequencies of the system are arbitrary large, and where resonances phenomena are known to occur for some values of the step size. For example, considering the case of the Schr\"odinger equation on the one dimensional torus, the eigenvalues of $H_0^{(K)}$ range from $1$ to $K^2$ and the assumption $h K^{2} << 1$ used in the finite dimensional situation becomes drastically restrictive in practice.  

Recently, many progresses have been made in the understanding of the long time behaviour of numerical methods applied to Hamiltonian PDEs. A first result using normal form techniques was given by {\sc Dujardin \& Faou} in \cite{DF07} for the case of the linear Schr\"odinger equation with small potential. Concerning the nonlinear case, results exists by {\sc Cohen, Hairer \& Lubich}, see \cite{CHL08b,CHL08c}, for the wave equation and {\sc Gauckler \& Lubich}, see \cite{GL08a,GL08b}, for the nonlinear Schr\"odinger equation using the technique of modulated Fourier expansion. However to be valid these results use non-resonance conditions that are generically satisfied only under CFL conditions linking the step-size $h$ and the highest frequencies of the discretized Hamiltonian PDE.

Normal form techniques have proven to be one of the most important tools for the understanding of the long time behaviour of Hamiltonian PDE (see \cite{Bam03, BG06, Greb07, Bam07, BDGS, GIP}). Roughly speaking, the dynamical consequences of such results are the following: starting with an small initial value of size $\varepsilon$ in a Sobolev space $H^s$, then the solution remains small in the same norm over long time, namely for time $t \leq C_r \varepsilon^{-r}$ for arbitrary $r$ (with a constant $C_r$ depending on $r$). Such results hold under {\em generic} non resonance conditions on the frequencies of the underlying linear operator $H_0$ associated with the Hamiltonian PDE, that are valid in a wide number of  situations (nonlinear Schr\"odinger equation on a torus of dimension $d$ or with Dirichlet boundary conditions, nonlinear wave equation with periodic or Dirichlet conditions in dimension 1, Klein Gordon equation on spheres or Zoll manifolds.).

This work is the first of a series of two.

In this paper, we consider full discretizations of the Hamiltonian PDE, with a spectral discretization parameter $K$ that is {\em finite but large}.  We show that under the hypothesis $K \leq \varepsilon^{-\sigma}$ for some constant $\sigma$ depending on the precision degree $r$ then the {\em actions} of the initial value are almost preserved over a very large number of iterations $n \leq C_r \varepsilon^{-r}$, provided the initial solution is small (of order $\varepsilon$) in $L^2$ norm. These actions can be interpreted as the oscillatory energies corresponding to an eigenvalue of $H_0^{(K)}$. Moreover, the $L^2$ norm of this numerical solution remains small for this large number of iterations. 

The method used in this situation is by essence a finite dimensional Birkhoff normal form result (explaining why we work here essentially with the $L^2$ norm). Using a {\em generic} non resonance condition on the step size that turns out to be valid for many equations and that is independent on $K$, we mainly show that we can take  $K$ asymptotically large without altering the nature of the classical result. 
Our main result is given by Theorem \ref{Tmain}. 

Roughly speaking, the method consists in applying techniques that are now standard in normal form theory, by tracking the dependence in $K$ of the constants appearing in the estimates. The use of a non resonance condition that is independent of $K$ is however crucial, and reflects the infinite dimensional nature of the initial continuous problem without space approximation. 

In some sense, the second paper \cite{FGP2} studies the case where $K > \varepsilon^{-\sigma}$ by considering the splitting method where no discretization in space is made (i.e. $K = +\infty$). The techniques used involve the abstract framework developed in \cite{BG06, Greb07, Bam07}. However, instead of being valid for the (exact) abstract splitting \eqref{E001}, we have to consider {\em rounded} splitting methods of the form
\begin{equation}
\label{E002}
\Pi_{\eta,s}\circ \varphi_{H_0}^h\circ \varphi_P^h
\end{equation}
where $\Pi_{\eta,s}$ puts to zero all the frequencies $\xi_j$ whose weighted energy $|j|^{2s}|\xi_j|^{2s}$ in the Sobolev space $H^s$ is smaller than a given threshold $\eta^2$. Hence, for small $\eta$, \eqref{E002} is very close to the exact splitting method \eqref{E001}. The good news is that this threshold can be taken of the order $\varepsilon^r$, making this projection $\Pi_{\eta,s}$ very close to the identity, and in any case producing an error that is far beyond the round-off error in a computer simulation (particularly for large $s$).

\section{Description of the method}

Before going on into the precise statements and proofs of this work, we would like to give  tentative explanations of the restrictions observed in comparison with the continuous case. 

The method used in \cite{BG06} to prove the long-time conservation of Sobolev norms and the associated weighted actions for small data is to start from a Hamiltonian $H = H_0 + P$ depending on an infinite number of variable $(\xi_j,\eta_k)$, $j, k \in \N$, and for a fixed number $r$, to construct a Hamiltonian transformation $\tau$ close to the identity, and such that in the new variable, the Hamiltonian can be written
\begin{equation}
\label{EnfBG}
H_0  + Z + R
\end{equation}
where $Z$ is a real Hamiltonian depending only on the action $I_j = \xi_j\eta_j$ and $R$ a real Hamiltonian having a zero of order $r$. 

The key for this construction is an induction process with, at each step, the resolution of an homological equation of the form 
\begin{equation}
\label{HomBG}
\{ H_0, \chi\} + Z = G
\end{equation}
where $G$ is a given homogeneous polynomial of order $n$, and where $Z$, depending only on the actions, and $\chi$ are unknown. Assume that $G$ is of the form 
$$
G = G_{\jb\kb} \,\xi_{j_1}\cdots \xi_{j_p}\eta_{k_1}\cdots \eta_{k_q}
$$
where $G_{\jb\kb}$ is a coefficient, $\jb = (j_1,\ldots,j_p) \in \N^p$ and $\kb = (k_1,\ldots,k_q) \in \N^q$. 
Then it is easy to see that the equation \eqref{HomBG} can be written 
\begin{equation}
\label{HomBG2}
\Omega(\jb,\kb) \chi_{\jb\kb} + Z_{\jb\kb} = G_{\jb\kb}
\end{equation}
where 
$$
\Omega(\jb,\kb) = \omega_{j_1} + \cdots + \omega_{j_p} - \omega_{k_1} -\cdots - \omega_{j_q}
$$
and where $Z_{\jb\kb}$ and $\chi_{\jb\kb}$ are unknown coefficients. 

It is clear that for $\jb = \kb$ (up to a permutation), we have $\Omega(\jb,\kb) = 0$ which imposes $Z_{\jb\kb} = G_{\jb\kb}$. When $\jb \neq \kb$ (taking into account the permutation), the solution of \eqref{HomBG2} relies on a non resonance conditions on the small divisors $\Omega(\jb,\kb)^{-1}$. 

In \cite{BG06}, {\sc Bambusi \& Gr\'ebert} use a non resonance condition of the form 
\begin{equation}
\label{EnonresBG}
\forall\, \jb \neq \kb,\quad | \Omega(\jb,\kb) | \geq \gamma \mu(\jb,\kb)^{-\alpha}
\end{equation}
where $\mu(\jb,\kb)$ denotes the third largest integer amongst $|j_1|,\ldots,|k_q|$. They moreover show that such a condition is guaranteed in a large number of situations (see \cite{BG06}, \cite{Greb07} or \cite{Bam07} for precise results). 

 Considering now the splitting method $\varphi_{H_0}^h \circ \varphi_{P}^h$, we see that we cannot work directly at the level of the Hamiltonian. To avoid this difficulty, we embed the splitting into the family of applications
$$
[0,1] \ni \lambda\mapsto \varphi_{H_0}^h \circ \varphi_{hP}^\lambda
$$
and we derive this expression with respect to $\lambda$, in order to work in the tangent space, where it is much easier to identify real Hamiltonian than unitary flows. 

This explains why we deal here with time-dependent Hamiltonian. Note that we do not expand the operator $\varphi_{H_0}^h$  in powers of $h$, as this would yields positive powers of the unbounded operator $H_0$ appearing in the series. Unless a CFL condition is employed, this methods do not give the desired results (and do not explain the resonance effects observed for some specific values of $h$). 

Now, instead of \eqref{HomBG}, the Homological equation appearing for the splitting methods is given in a discrete form
\begin{equation}
\label{HomFGP}
\chi \circ \varphi_{H_0}^h - \chi + Z = G. 
\end{equation}
In terms of coefficients, this equations yields
$$
(e^{ih\Omega(\jb,\kb)} - 1) \chi_{\jb\kb} + Z_{\jb\kb} = G_{\jb\kb}.
$$
The main difference with \eqref{HomBG2} is that we have to avoid not only the indices $(\jb,\kb)$ so that $\Omega(\jb,\kb) = 0$, but all of those for which $h \Omega(\jb,\kb) = 2 m \pi$ for some (unbounded) integer $m$.

In the case of a fully discretized system for which $\nabla_{z_j} P \equiv 0$ for $|j| > K$, then under the CFL-like condition of the form $h K^m\leq C$ where $m$ depends on the growth of the eigenvalues of $H_0$ and $C$ depends on $r$, then we have $|h \Omega(\jb,\kb)| \leq  \pi$, and hence 
\begin{equation}
\label{EnonresGL}
|e^{ih\Omega(\jb,\kb)} - 1| \geq h \gamma \mu(\jb,\kb)^{-\alpha}
\end{equation}
\eqref{EnonresGL}
is then a consequence of \eqref{EnonresBG}. Under this assumption, we can apply the same techniques used in  \cite{BG06} and draw the same conclusions. This is the kind of assumption made in \cite{CHL08b} and \cite{GL08b}.

The problem with \eqref{EnonresGL} is that it is non generic in $h$ outside the CFL regime. For example, in the case of the Schr\"odinger equation, the frequencies of the operator $H_0$ are such that  $\omega_j \simeq j^2$. Hence, for large $N$, if $(j_1,\ldots,j_p,k_1,\ldots,k_q)$ is such that $j_1 = N +1$, $k_1 = N$ and all the other are of order $1$ ($N$ is large here), we have $\Omega(\jb,\kb) \simeq  (N+1)^2 - N^2 \simeq 2N$. Hence, 
$$
|e^{ih\Omega(\jb,\kb)} - 1| \simeq | e^{2ihN} - 1 | 
$$
cannot be assumed to be greater than $h \gamma \mu(\jb,\kb)^{-\alpha} \simeq h$ for all (large) $N$. Note that a generic hypothesis on $h$ would be here that this small divisor is greater than $h \gamma N^{-\alpha}$ for some constants $\gamma$ and $\alpha$. This example shows that we cannot control the small divisors $|e^{ih\Omega(\jb,\kb)} - 1|$ associated with the splitting scheme by the {\em third largest} integer in the multi index (which is actually of order 1 in this case), but by the {\em largest}. 

Using a generic condition on $h \leq h_0$, we prove in \cite{FGP2} a normal form result and show that the flow is conjugated to the flow of a Hamiltonian vector field of the form \eqref{EnfBG}, but where $Z$ now contains terms depending only on the actions, and supplementary terms containing  at least two large indices. Here, large means greater than $\varepsilon^{-\sigma}$ where $\sigma$ depends on $r$. 

In the case of a full discretization of the Hamiltonian PDE with a spectral discretization parameter $K$, we thus see that if $K \leq \varepsilon^{-\sigma}$ then the normal form term $Z$ actually depends only on the actions, as the high frequencies greater that $\varepsilon^{-\sigma}$ are not present. This is essentially the result of this paper. 

In the case where $K > \varepsilon^{-\sigma}$, the normal form result that we obtain can be interpreted as follows: the non conservation of the actions can only come from two high modes (of order greater than $\varepsilon^{-\sigma}$) interacting together and contaminating the whole spectrum. The role of the projection operator $\Pi_{\eta,s}$ is to destroy these high modes at each step. As we can take $\eta = \varepsilon^r$, the error induced is very small, and in particular, far beyond the round-off error in the numerical simulation. This is mainly the result in \cite{FGP2}.

The differences between the present work and \cite{FGP2} lie in the techniques involved: In this work, the system considered are {\em large} but {\em finite} dimensional systems, and all the hypothesis made on the nonlinearity can be expressed using elementary conditions similar to those used in the finite dimensional case. 
In \cite{FGP2}, we study $K = +\infty$, which requires much more elaborate technical tools in the spirit of \cite{BG06, Greb07, Bam07}.

\section{Setting of the problem}

\subsection{Hamiltonian formalism}

We set $\Nc = \Z^d$ or $\N^d$. For $a = (a_1,\ldots,a_d) \in \Nc$, we set
$$
|a| = \max_{i = 1,\ldots,d} |a_i|. 
$$
Let $K \in \N$, and let $\Nc_K$ a finite subset of $\Nc$, included in the ball $\{ a \in \Nc\, | \, |a| \leq K\,\}$. 
Typically, we can take $\Nc_K$ of the form  $[-K,\ldots,K]^d \subset \Z^d$ or  $[0,\ldots,K]^d \subset \N^d$ or a sparse set of the form  (see for instance \cite{G07,L08})
$$
\Nc_K= \{ \, a = (a_1,\ldots,a_d)  \in \Z^d\; | \;  ( 1 + |a_1|)\cdots(1 + |a_d|) \leq K\, \} \subset \Z^d.
$$

We consider the set of variables $(\xi_a,\eta_b) \in \C^{\Nc_K} \times \C^{\Nc_K}$ equipped with the symplectic structure
\begin{equation}
\label{Esymp}
i \sum_{a \in \Nc_K} \dd \xi_a \wedge \dd \eta_a. 
\end{equation}
We define the set $\Zc_K = \Nc_K \times \{ \pm 1\}$. For $j = (a,\delta) \in \Zc_K$, we define $|j| = |a|$ and we denote by $\overline{j}$ the index $(a,-\delta)$. 

We then define the variables
$(z_j)_{j \in \Zc_K} \in \C^{\Zc_K}$ by the formula
$$
j = (a,\delta) \in \Zc_K  \Longrightarrow 
\left\{
\begin{array}{rcll}
z_{j} &=& \xi_{a}& \mbox{if}\quad \delta = 1,\\[1ex]
z_j &=& \eta_a & \mbox{if}\quad \delta = - 1,
\end{array}
\right.
$$
By abuse of notation, we often write $z = (\xi,\eta)$ to denote such an element. 

We set 
$$
\Norm{z}{}^2 := \sum_{j \in \Zc_K} |z_j|^2
$$
and for any $\rho >0$, 
$$
B_K(\rho) = \{ \,  z \in \C^{\Zc_K}\, | \, \Norm{z}{} \leq \rho\,\}.
$$
Note that in the case where $K = +\infty$, we set by convention $\Zc_K = \Zc = \Nc \times \{\pm 1 \}$ and the previous norm defines a Hilbert structure on  $\ell^2_\Zc$. We denote by 
$$
\Pi_K: \ell^2_\Zc \to \big(\C^{\Zc_K},\Norm{\cdot}{}\big)
$$
the natural projection. 

Let $\Uc_K$ be a an open set of $\C^{\Zc_K}$. For a function $F$ in $\mathcal{C}^1(\Uc_K,\C)$, we define its gradient as 
$$
\nabla F(z) = \left( \frac{\partial F}{\partial z_j}\right)_{j \in \Zc_K}
$$
where by definition, we set for $j = (a,\delta) \in \Nc_K \times \{ \pm 1\}$, 
$$
 \frac{\partial F}{\partial z_j} =
  \left\{\begin{array}{rll}
 \displaystyle  \frac{\partial F}{\partial \xi_a} & \mbox{if}\quad\delta = 1,\\[2ex]
 \displaystyle \frac{\partial F}{\partial \eta_a} & \mbox{if}\quad\delta = - 1.
 \end{array}
 \right.
$$
Let $H(z)$ be a function defined on $\Uc_K$. If $H$ is smooth enough, we can associate with this function the Hamiltonian vector field $X_H(z)$ defined as
$$
X_H(z) = J \nabla H(z) 
$$
where $J$ is the symplectic operator induced by the symplectic form \eqref{Esymp}.

For two functions $F$ and $G$, the Poisson Bracket is defined as
$$
\{F,G\} = \nabla F^T J \nabla G = i \sum_{a \in \Nc_K} \frac{\partial F}{\partial \eta_j}\frac{\partial G}{\partial \xi_j} -  \frac{\partial F}{\partial \xi_j}\frac{\partial G}{\partial \eta_j}.  
$$

We say that $z\in \C^{\Zc_K}$ is {\em real} when $z_{\overline{j}} = \overline{z_j}$ for any $j\in \Zc_K$. In this case, $z=(\xi,\bar\xi)$ for some $\xi_K\in \C^{\Nc_K}$. Further we say that a Hamiltonian function $H$ is 
 {\em real } if $H(z)$ is real for all real $z$. 

With a given function $H \in  \mathcal{C}^{\infty}(\Uc_K,\C)$, we associate the Hamiltonian system
$$
\dot z = J \nabla H(z)
$$
which can be written
\begin{equation}
\label{Eham2}
\left\{
\begin{array}{rcll}
\dot\xi_a &=& \displaystyle - i \frac{\partial H}{\partial \eta_a}(\xi,\eta) & a \in \Nc_K\\[2ex]
\dot\eta_a &=& \displaystyle i \frac{\partial H}{\partial \xi_a}(\xi,\eta)& a \in \Nc_K.  
\end{array}
\right.
\end{equation}
In this situation, we define the flow $\varphi_H^t(z)$ associated with the previous system (for times $t \geq 0$ depending on $z \in \Uc_K$). Note that if $z = (\xi,\bar \xi)$ and $H$ is real, the flow $(\xi^t,\eta^t) = \varphi_H^t(z)$, for all time where it is defined, satisfies the relation $\xi^t = \bar {\eta}^t$, where $\xi^t$ is solution of the equation 
\begin{equation}
\label{Eham1}
\dot\xi_a = - i \frac{\partial H}{\partial \eta_a}(\xi,\bar\xi), \quad a \in \Nc_K. 
\end{equation}
In this situation, introducing the real variables $p_a$ and $q_a$ such that
$$
\xi_a = \frac{1}{\sqrt{2}} (p_a + i q_a)\quad \mbox{and}\quad \bar{\xi}_a =  \frac{1}{\sqrt{2}} (p_a - i q_a),
$$
the system \eqref{Eham1} is equivalent to the system
$$
\left\{
\begin{array}{rcll}
\dot p_a &=& \displaystyle -  \frac{\partial \tilde{H}}{\partial q_a}(q,p) & a \in \Nc_K\\[2ex]
\dot q_a &=& \displaystyle  \frac{\partial \tilde{H}}{\partial p_a}(q,p),& 	a \in \Nc_K.  
\end{array}
\right.
$$
where $\tilde{H}(q,p) = H(\xi,\bar\xi)$. 

Note that the flow $\tau^t = \varphi_\chi^t$ of a real Hamiltonian $\chi$ defines a symplectic map, i.e.  satisfies for all time $t$ and all point $z$ where it is defined
\begin{equation}
\label{Esympl}
(D_z  \tau^t)_z^T J (D_z\tau^t)_z = J
\end{equation}
where $D_z$ denotes the derivative with respect to the initial conditions. 

The following result is classical: 
\begin{lemma}
\label{Lchange}
Let $\Uc_K$ and $\Wc_K$ be two domains of $\C^{\Zc_K}$, and let $\tau = \varphi_\chi^1 \in \mathcal{C}^{\infty}(\Uc_K,\Wc_K)$ be the flow of the real smooth Hamiltonian $\chi$. 
Then for $H \in \mathcal{C}^{\infty}(\Wc_K,\C)$, we have 
$$
\forall\, z \in \Uc\quad
X_{H \circ \tau}(z) = 
(D_z\tau(z))^{-1}X_H(\tau(z)).
$$
Moreover, if $H$ is a real Hamiltonian, $H \circ \tau$ is a real Hamiltonian. 
\end{lemma}

\subsection{Hypothesis}

We describe now the hypothesis needed on the Hamiltonian $H$. 

In the following, we consider an infinite set of frequencies $(\omega_a)_{a \in \Nc}$ satisfying
\begin{equation}
\label{Eboundomega}
\forall\, a \in \Nc, \quad |\omega_a| \leq C |a|^m
\end{equation}
for some constants $C > 0$ and $m > 0$.

Let $\Uc$ be an open domain of $\ell^2(\C^{\Zc})$ containing the origin, and let $\Uc_K = \Pi_K \Uc$ its  projection onto $\C^{\Zc_K}$. 

We consider the collection of Hamiltonian functions 
\begin{equation}
\label{Edecomp}
H^{(K)} = H_0^{(K)} + P^{(K)},\quad K \geq 0, 
\end{equation}
with
$$
H_0^{(K)} = \sum_{a \in \Nc_K} \omega_a I_a(z)
$$
where for all $a\in \Nc_K$, 
\begin{equation}
\label{Eaction}
I_a(z) = \xi_a \eta_a
\end{equation}
are the {\em actions} associated with $a\in \Nc_K$. Note that if $z = (\xi,\bar\xi)$, then $I_a(z) = |\xi_a|^2$. 

We moreover assume that the functions $P^{(K)} \in \mathcal{C}^{\infty}(\Uc_K,\C)$ are {\em real}, of {\em order at least 3},  and  satisfy the following: For all $\ell > 1$, there exists constants $C(\ell) \geq 0$ and $\beta(\ell)\geq 0$ such that for all $K \geq 1$, $(j_1,\cdots,j_\ell) \in \Zc_K^\ell$ and $z \in \Uc_K$, the following estimate holds:
\begin{equation}
\label{EestP}
\left|
\frac{\partial P^{(K)}}{\partial z_{j_1} \cdots \partial z_{j_\ell}}(z)
\right|
\leq C(\ell) K^{\beta(\ell)}. 
\end{equation}

The Hamiltonian system \eqref{Eham2} can hence be written
\begin{equation}
\label{Eham3}
\left\{
\begin{array}{rcll}
\dot\xi_a &=& \displaystyle - i \omega_a \xi_a - i \frac{\partial P^{(K)}}{\partial \eta_a}(\xi,\eta) & a \in \Nc_K\\[2ex]
\dot\eta_a &=& \displaystyle i \omega_a \eta_a + i \frac{\partial P^{(K)}}{\partial \xi_a}(\xi,\eta)& a \in \Nc_K.  
\end{array}
\right.
\end{equation}
Denoting by $\varphi_Q^t$ the exact flow of a Hamiltonian flow, splitting methods are based on the approximation
$$
\varphi_{H^{(K)}}^h \simeq \varphi_{H_0^{(K)}}^h \circ \varphi_{P^{(K)}}^h
$$ 
for a small time step $h > 0$. Note that in this case, the exact flow of $H_0^{(K)}$ is explicit and given by 
$$
\varphi_{H_0^{(K)}}^h(\xi,\eta) = (e^{-i\omega_a h} \xi_a, e^{i \omega_a h} \eta_a)_{a \in \Nc_K}
$$
while the calculation of $\varphi_{P^{(K)}}^h$ requires the solution of an ordinary differential equation, whose solution is often given explicitely (see the examples below). 

The goal of this paper is the study of the long-time behavior of the numerical solution $z^n$ given by \eqref{Ezn} for large number $n$ of iterations. 

\begin{remark}
Note that no hypothesis is made here concerning the preservation of the $L^2$ norm by the flow of \eqref{Eham3}. 
\end{remark}
\subsection{Non resonance condition}
\label{SSak}

In the following, for $\jb = (j_1,\ldots,j_r) \in \Zc_K^r$ with $r \geq 1$, we use the notation 
$$
z_\jb = z_{j_1}\cdots z_{j_r}. 
$$
Moreover, for $\jb = (j_1,\ldots,j_r) \in \Zc_K^r$ with $j_i = (a_i,\delta_i) \in \Nc_K \times\{ \pm 1\}$ for $i = 1,\ldots,r$, we set
$$
\overline \jb = (\overline{j}_1,\ldots,\overline j_r)\quad\mbox{with}\quad \overline{j}_i = (a_i,-\delta_i), \quad i = 1,\ldots,r,
$$
and we define 
$$
\Omega(\jb) = 
\delta_1\omega_{a_1} + \cdots  + \delta_r\omega_{a_r}. 
$$
We say that $\jb \in \Zc_K^r$ depends only of the action and we write $\jb \in \Ac_K^r$ if $r$ is even and if we can write (up to a permutation of the indexes)
$$
\forall\, i = 1,\ldots r/2,\quad
j_{i} = (a_i,1), \quad\mbox{and}\quad j_{i + r/2} = (a_i,-1)
$$
for some $a_i \in \Nc_K$.   
Note that in this situation, 
$$
\begin{array}{rcl}
z_\jb = z_{j_1}\cdots z_{j_r} &=& \xi_{a_1}\eta_{a_1} \cdots \xi_{a_{r/2}} \eta_{a_{r/2}}\\[2ex]
&=& I_{a_1}(z) \cdots I_{a_{r/2}}(z)
\end{array}
$$
where for all $a \in \Nc_K$, $I_a(z)$ denote the actions associated with $a$ (see \eqref{Eaction}). 
For odd $r$, $\Ac_r$ is the empty set.

We will assume now that the step size $h$ satisfies the following property: 
\begin{hypothesis}\label{H1}
For all $r \in \N$, there exist constants $\gamma^*$ and $\alpha^*$ such that for all $K \in \N^*$,
\begin{equation}
\label{nonres2}
(j_1,\ldots,j_r) \in \Zc_K^r \backslash \Ac_K^r\quad \Longrightarrow \quad |1 - e^{ih\Omega(\jb)}| \geq \frac{h \gamma^*}{K^{\alpha^*}}. 
\end{equation}
\end{hypothesis}

The following Lemma \ref{Lnonres} shows that condition \eqref{nonres2} is {\em generic} in the sense that it is satisfied for a large set of $h \leq h_0$ (and in particular independently of $K$), provided that the frequencies $\omega_a$ satisfy a non resonance condition that we state now (see \cite{HLW, Shan00} for similar statements):

\begin{hypothesis}\label{H2}
For all $r \in \N$, there exist  constants $\gamma(r)$ and $\alpha(r)$ such that $\forall\, K \in \N^*$, 
\begin{equation}
\label{nonres}
(j_1,\ldots,j_r) \in \Zc_K^r \backslash \Ac_K^r \quad   \Longrightarrow |\Omega(\jb)| \geq \frac{ \gamma}{K^{\alpha}}. 
\end{equation}
\end{hypothesis}

In the next section, we will check that condition in different concrete cases. 

\begin{lemma}
\label{Lnonres}
Assume that Hypothesis \ref{H2} holds, and let $h_0$ and $r$ be given numbers. Let $\gamma$ and $\alpha$ be such that \eqref{nonres} holds and assume that $\gamma^* \leq  (2/\pi) \gamma$, $\alpha^* \geq \alpha + m\sigma + r$ with $\sigma> 1$ and $m$ the constant appearing in \eqref{Eboundomega}, then we have 
$$
\mbox{\rm meas}\{\, h < h_0\, | \, h \mbox{ does not satisfy }  \eqref{nonres2}\, \} \leq C \frac{\gamma^*}{\gamma} h_0^{1+\sigma}
$$
where $C$ depends on $\sigma$ and  $r$. As a consequence the set  
$$
Z(h_0) = \{\, h < h_0\, | \, h \mbox{ satisfies Hypothesis }  \ref{H1}\, \}
$$
is a dense open subset of $(0,h_0)$.

\end{lemma}

The proof of this lemma is given in \cite[Lemma 4.6]{FGP2}.

\section{Statement of the result and applications}

\subsection{Main results}

\begin{theorem}
\label{Tmain}
Assume that $P^{(K)}$ and $h < h_0$ satisfy the previous hypothesis. Let $r \in \N^*$  be fixed. There exist constants $\sigma$, $C$ and $\varepsilon_0$ depending only on  $r$, $h_0$ and the constants $\beta(\ell)$ and $C(\ell)$, $\ell = 0,\ldots r$ in \eqref{EestP},
such that 
the following holds: For all $\varepsilon < \varepsilon_0$ and $K \leq \varepsilon^{-\sigma}$, and for all $z^0$ real such that
$$
\Norm{z^0}{} \leq \varepsilon 
$$ 
if we define 
\begin{equation}
\label{Eseuil}
z^n = \big(\varphi_{H_0^{(K)}}^h \circ \varphi_{P^{(K)}}^h \big)^n (z^0)
\end{equation}
then
for all $n$,  $z^n$ is still real, 
and moreover
\begin{equation}
\label{Eresnorm}
\Norm{z^n}{} \leq 2 \varepsilon\quad \mbox{for}\quad n \leq \frac{1}{\varepsilon^{r-1}},
\end{equation}
and
\begin{equation}
\label{Eresact}
\forall\, a \in \Nc_K, \quad | I_a(z^n) - I_a(z^0)| \leq C\varepsilon^{5/2}\quad \mbox{for}\quad n \leq \frac{1}{\varepsilon^{r-2}}
\end{equation}
\end{theorem}

The proof of this result relies on the following Birkhoff normal form result, whose proof is postponed to Section \ref{SProof}: 
\begin{theorem}
\label{TNF}
Assume that that $P^{(K)}$ and $h < h_0$ satisfy hypothesis \eqref{EestP} and \eqref{nonres2}. Let $r \in \N^*$ be fixed. Then there exists
constants $\beta$ and $C$ depending on $r$, $h_0$, $\beta(\ell)$ and $C(\ell)$, $\ell = 0,\ldots r$ in \eqref{EestP} and a canonical transformation $\tau_K$ from $B_K(\rho)$ into $B_K(2\rho)$ with $\rho = (CK)^{-\beta}$ satisfying for all $z \in B_K(\rho)$, 
\begin{equation}
\label{Eesttau}
\Norm{\tau_K(z) - z}{} \leq (CK)^\beta \Norm{z}{}^2
\quad \mbox{and}\quad 
\Norm{\tau_K^{-1}(z) - z}{} \leq (CK)^\beta \Norm{z}{}^2, 
\end{equation}
satisfying the following result: For all $z \in B_K(\rho)$, 
$$
\tau_{K}^{-1} \circ \varphi_{H_0^{(k)}}^h \circ \varphi_{P^{(K)}}^h \circ \tau_K(z) = \varphi_{H_0^{(K)}}^h \circ \psi_K (z)
$$
where $\psi_K$ satisfies: 
\begin{itemize}
\item $\psi_K(z)$ is real if $z$ is real, 
\item For all $z \in B_K(\rho)$, 
\begin{equation}
\label{Enorme}
\Norm{\psi_K(z) - z}{}  \leq (CK)^\beta \Norm{z}{}^{r}. 
\end{equation}
\item For all $z \in B_K(\rho)$, 
\begin{equation}
\label{Eactions}
|I_a(\psi_K(z)) - I_a(z)|  \leq (CK)^\beta \Norm{z}{}^{r+1}. 
\end{equation}
\end{itemize}
\end{theorem}

\begin{Proofof}{Theorem \ref{Tmain}}
First, let us note that as the Hamiltonian functions $H_0^{(K)}$ and $P^{(K)}$ are real Hamiltonians, it is clear that there exist $\xi^n \in \C^\Nc$ such that for all $n$, we have $z^n = (\xi^n,\bar\xi^n)$, that is $z^n$ is real. 

Let $\beta$ given by Theorem \ref{TNF} and let $\sigma = 1/(2\beta)$. 
We have for $K \leq \varepsilon^{-\sigma}$, 
$$
(C K) ^{\beta} \leq C^\beta \varepsilon^{-1/2}. 
$$

Let $\tau_K$ be defined by Theorem \ref{TNF}, and let $y^n = \tau_K^{-1}(z^n)$. 
Using the property of $\tau_K$, we see that $y^n$ is real, i.e. we have $y^n = (\zeta^n,\bar\zeta^n)$ for all $n$.  By definition, we have
\begin{equation}\label{star}
\forall\, n \geq 0,\quad 
y^{n+1} =  \big(\varphi_{H_0^{(K)}}^h \circ \psi_K\big) (y^n). 
\end{equation}
Using the fact that $K \leq \varepsilon^{-\sigma}$ and \eqref{Eesttau}, the transformation $\tau_K$ in the previous Theorem satisfies the following: For all $z$ such that $\Norm{z}{} \leq 2\varepsilon$, 
\begin{equation}
\label{Etransfo}
\Norm{\tau_K^{-1}(z) - z}{} \leq   C^{\beta}\varepsilon^{-1/2} \Norm{z}{}^2 
\leq 4 C^\beta \varepsilon^{3/2} \leq \textstyle\frac14 \varepsilon
\end{equation}
provided $\varepsilon_0$ is sufficiently small. 
Hence we have 
$\Norm{y^0}{} = \Norm{\tau_K^{-1}(z^0)}{} \leq \frac54 \varepsilon$. 

Note that we have $\rho = (CK)^{-\beta} \geq C^{-\beta}\varepsilon^{1/2} \geq 2 \varepsilon$ provided that $\varepsilon_0$ is small enough. 
Using \eqref{Enorme} we get that as long as $\Norm{y^n}{} \leq 2 \varepsilon$, we have
$$
\Norm{y^{n+1}}{} \leq \Norm{y^n}{} + (CK)^{\beta} \Norm{y^n}{}^{r} \leq \Norm{y^n}{} + 2^r C^\beta \varepsilon^{r -1/2} 
$$
By induction, we thus see that for 
$$
n \leq 2^{-r-1} C^{-\beta} \varepsilon^{3/2 - r}
$$
we have $\Norm{y^{n}}{} \leq \frac74 \varepsilon \leq 2 \varepsilon$. Assuming that $\varepsilon_0$ is such that $2^{-r-1} C^{-\beta}\varepsilon_0^{1/2} \leq 1$, this shows that for $n \leq \varepsilon^{1-r}$ we have  $\Norm{y^{n}}{} \leq \frac74 \varepsilon$. 
Using \eqref{Eesttau} and an inequality similar to \eqref{Etransfo}, we conclude that 
$$
\Norm{z^{n}}{} \leq 2\varepsilon,\quad \mbox{for}\quad n \leq \frac{1}{\varepsilon^{r-1}}
$$
which yields to \eqref{Eresnorm}. 

Now using \eqref{Eactions} and the fact that $\Norm{y^n}{} \leq 2 \varepsilon$ we see that 
for $n \leq \varepsilon^{1-r}$ we have 
$$
\forall\, a \in \Nc_K, \quad | I_a(y^{n+1}) - I_a(y^n)| \leq 2^{r+1} C^\beta \varepsilon^{r+1/2}
$$
whence 
$$
\forall\, a \in \Nc_K, \quad | I_a(y^{n}) - I_a(y^0)| \leq 2^{r+1} C^\beta n \varepsilon^{r+1/2}
$$
Now we have for all $a \in \Nc_K$
$$
| I_{a}(y^n) - I_a(z^n)| = \big| |\zeta_a^n|^2 - |\xi_a^n|^2 \big| = \big| |\zeta_a^n| - |\xi_a^n| \big| \times  \big| |\zeta_a^n| + |\xi_a^n| \big|, 
$$
whence
$$
| I_{a}(y^n) - I_a(z^n)| \leq |\zeta_a^n - \xi_a^n| (\Norm{y^n}{} + \Norm{z^n}{}) \leq \Norm{\tau_K(y^{n}) - y^n}{}(\Norm{y^n}{} + \Norm{z^n}{}). 
$$
Using \eqref{Eesttau} we see that for all $n \leq \varepsilon^{1-r}$ and all $a \in \Nc_K$, 
$$
\Norm{\tau_K(y^{n}) - y^n}{} \leq  4 C^{\beta} \varepsilon^{3/2}.
$$
and hence, as $\Norm{z^n}{} \leq 2 \varepsilon$, 
$$
| I_{a}(y^n) - I_a(z^n)| \leq 8 C^{\beta} \varepsilon^{5/2}. 
$$
Using \eqref{Etransfo}, we thus see that 
$$
\forall\, n \leq \varepsilon^{1-r}\,\quad \forall\, a \in \Nc_K, \quad | I_a(z^{n}) - I_a(z^0)| \leq  2^{r + 4} C^{\beta} ( \varepsilon^{5/2} + n \varepsilon^{r+ 1/2})
$$
and this easily gives the result. 
\end{Proofof}

\subsection{Examples}
In this section we  present two examples, other examples like the Klein Gordon equation on the sphere (in the spirit of \cite{BDGS}) or the nonlinear Schr\"odinger operator with harmonic potential (in the spirit of \cite{GIP}) could also be considered with these technics.

\subsubsection{Schr\"odinger equation on the torus}

We first consider  nonlinear Schr\"odinger equations of the form
\begin{equation}
\label{E1}
i \partial_t \psi = - \Delta \psi + V \star\psi + \partial_2g(\psi,\bar \psi),\quad x \in \T^d
\end{equation}
where $V\in C^\infty(\T^d,\R)$,  $g\in C^\infty(\Uc,\C)$ where $\Uc$ is a neighborhood of the origin in $\C^2$. We assume that $g(u,\bar u) \in \R$, and that $g(u,\bar u) = \mathcal{O}(|u|^3)$. 
The corresponding Hamiltonian functional is given by 
$$
H(\psi,\bar\psi) = \int_{\T^d} | \nabla \psi | ^2 + \bar\psi (V \star \psi) + g(\psi,\bar\psi) \, \dd x. 
$$

Let $\phi_{a}(x) = e^{i a\cdot x}$, $a \in \Z^d$ be the Fourier basis on $L^2(\T^d)$. With the notation 
$$
\psi = \Big(\frac{1}{2\pi}\Big)^{d/2}\sum_{a\in \Z^d} \xi_{a} \phi_{a}(x) \quad \mbox{and}\quad
\bar \psi = \Big(\frac{1}{2\pi}\Big)^{d/2}\sum_{a\in \Z^d} \eta_{a}\bar \phi_{a}(x)
$$
 the (abstract) Hamiltonian associated with the equation \eqref{E1} can be formally written
\begin{equation}
\label{E2}
H(\xi,\eta) = 
\sum_{a \in \Z^d} \omega_a \xi_a \eta_a  + P(\xi,\eta). 
\end{equation}
Here  $\omega_{a} =|a|^2 +\hat V_a$ are the eigenvalues of the operator 
$$
\psi\mapsto  - \Delta\psi + V \star\psi,
$$
and we see that $\omega_{a}$ satisfy \eqref{Eboundomega} with $m = 2$.
Moreover, the nonlinearity function $P(\xi,\eta)$ posesses a zero of order $3$ at the origin. 
In this situation, it can be shown  that the Hypothesis \ref{H1} is fulfilled for a large set of potential $V$ (see \cite{BG06} or \cite{Greb07}). 

Following \cite{GL08a}, a space discretization of this equation using spectral collocation methods yields a problem of the form \eqref{Edecomp} with 
$$
\Nc_K = [-K,\ldots,K-1]^d
$$
and, with 
\begin{equation}
\label{Eschrinc}
u_K = \Big(\frac{1}{2\pi}\Big)^{d/2}\sum_{a\in \Nc_K} \xi_{a} \phi_{a}(x) \quad \mbox{and}\quad
v_K = \Big(\frac{1}{2\pi}\Big)^{d/2}\sum_{a\in \Nc_K} \eta_{a}\bar \phi_{a}(x)
\end{equation}
the nonlinearity reads
$$
P^{(K)}(\xi,\eta) = \int_{\T^d} \mathcal{Q}( g(u_K,v_K) ) \dd x
$$
where, for a function $\psi = (\frac{1}{2\pi})^{d/2}\sum_{a\in \Z^d} \psi_{a} \phi_{a}(x)$
$$
\mathcal{Q}( \psi ) = \sum_{a \in \Nc_K} \Big(\sum_{b \in \Z^d} \psi_{a + 2Kb} \Big)\phi_a(x)
$$
is the collocation operator associated with the points $x_a = \frac{\pi}{K}\in \T^d$, $a \in \Nc_K$. It is easy to verify that $P^{(K)}$ satisfies \eqref{EestP} for some constants $C(\ell)$ depending on $g$ and $\beta(\ell)$ depending on $g$ and the dimension $d$. 

Note that starting from a real initial value $u_K^0(x)$ (see \eqref{Eschrinc}) this system reduces to solving the system of ordinary differential equation
$$
\forall\, a \in \Nc_a\, \quad 
i \frac{\dd }{\dd t} u_K(x_a,t) = \mathcal{F}_{2K} \Omega \mathcal{F}_{2K}^{-1} u_K(x_a,t) + \partial_2 g (u_K(x_a,t), \overline{u_K(x_a,t)})
$$
where $\Omega$ is the matrix $(\omega_a)_{a \in \Nc_K}$ and $\mathcal{F}_{2K}$ the Fourier transform associated with $\Nc_K$. Note that in this case, the numerical solution \eqref{Ezn} is easily implemented: The linear part is diagonal and can be solved explicitely in the Fourier space, while the non-linear part is an ordinary differential equation with fixed parameter $x_a$ at each step. If moreover 
$g(u,\bar u) = G(|u|^2)$ for some real function $G$ then the solution of the nonlinear part is given explicitely by 
$\varphi_{P^{(K)}}^h(u) = \exp(-2ih G'(|u|^2) u$  using the fact that $|u|^2$ is constant for  a fixed point $x_a$. 

For high dimension $d$, the previous discretization is usually replaced by a discretization on sparse grid, i.e. with 
$$
\Nc_K= \{ \, a = (a_1,\ldots,a_d)  \in \Z^d\; | \; (1 + |a_1|)\cdots(1 +|a_d|) \leq K\, \} \subset \Z^d.
$$
As explained in \cite[Chap III.1]{L08}, methods exist to write the corresponding system under the symplectic form \eqref{Edecomp}, upon a possible loss in the approximation properties of the exact solution of \eqref{E1} by the solution of the discretized Hamiltonian $H^{(K)}$. Note that this does not influence the long time results proven here: In some sense we do not impose the nonlinearity $P^{(K)}(z)$ to approximate an exact nonlinearity $P(z)$. 

We give first a numerical illustration of resonance effects. We consider the equation
$$
i\partial_t \psi  = - \Delta \psi + V \star \psi + \varepsilon^2 |\psi|^2 \psi
$$
in the one dimensional torus $\T^1$, with initial value
$$
\psi_0(x) = \frac{2}{2 - \cos(x)}. 
$$
Note that this problem is equivalent to solving \eqref{E1} with a small initial value of order $\varepsilon$. 
We take $\varepsilon = 0.1$, $V$ with Fourier coefficients $\hat V_a = 2/(10 + 2 a^2)$ and $K = 200$ (i.e. $400$ collocation points). In Figure 1, we plot the actions of the numerical solution given by the Lie splitting algorithm \eqref{Ezn} in logarithmic scale. In the right we use the resonant stepsize $h = 2\pi/(\omega_7 - \omega_1) \simeq 0.17459\ldots$. In the left we plot the same result but with the non resonant stepsize $h = 0.174$. 
\begin{figure}[ht]
\label{figure1}
\begin{center}
\rotatebox{0}{\resizebox{!}{0.3\linewidth}{%
   \includegraphics{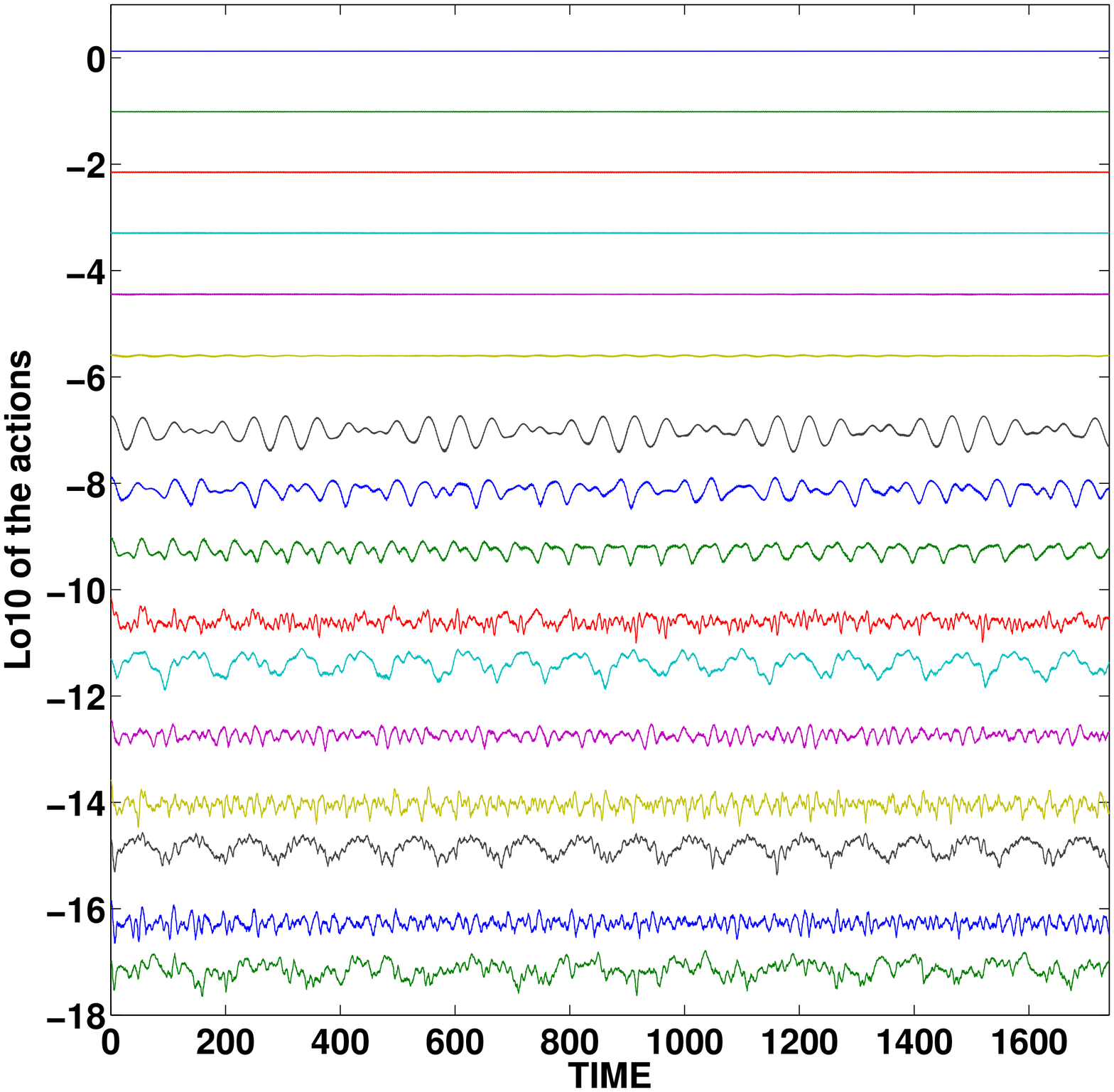}}} 
\rotatebox{0}{\resizebox{!}{0.3\linewidth}{%
   \includegraphics{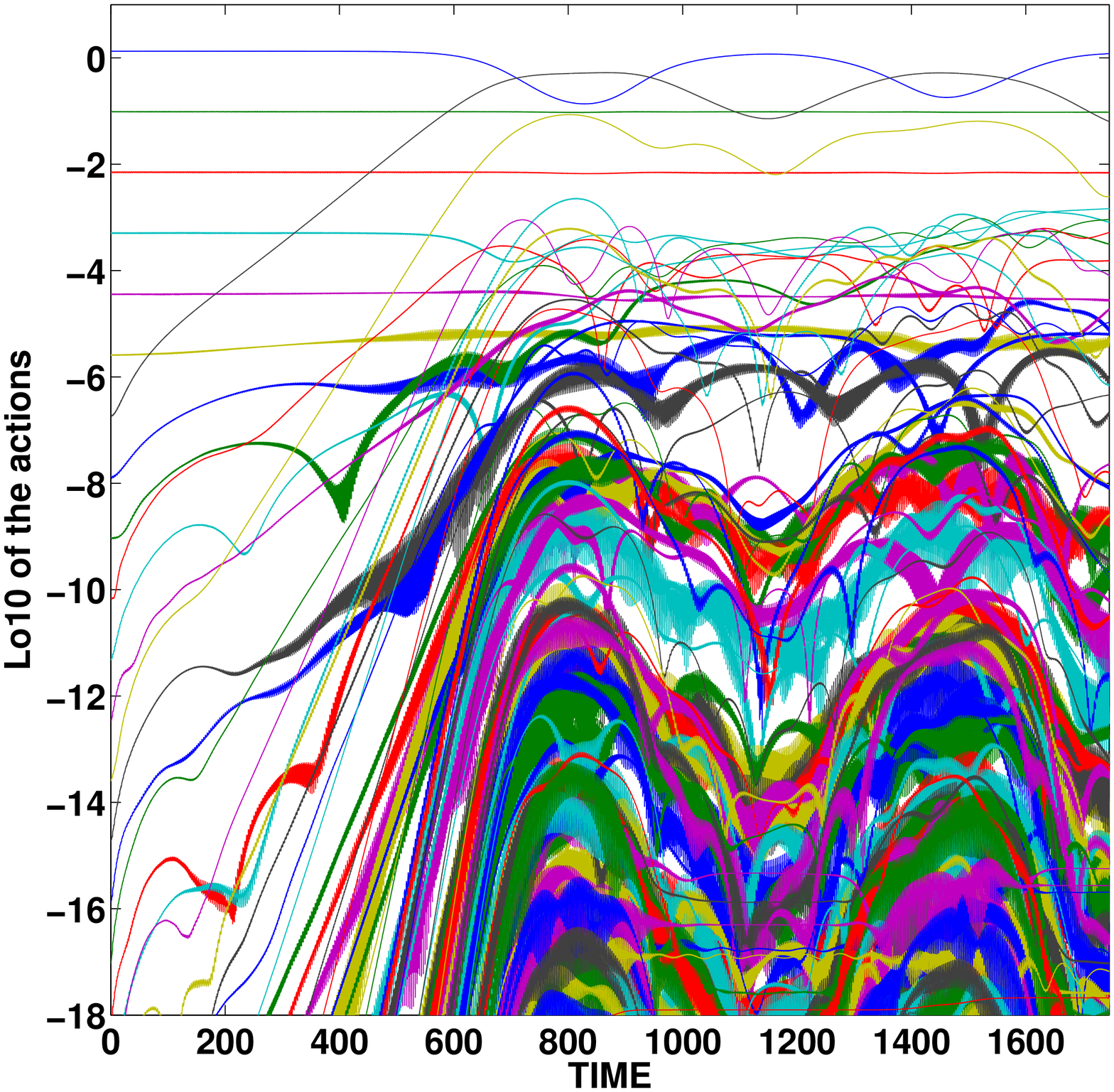}}} 
\caption{\small{Plot of the actions for non-resonant and resonant step size.} }
\end{center}
\end{figure}

In Figure 2, we show the long time almost conservation of the action in the case where $h = 0.1$ (non resonant), and $\varepsilon = 0.1$ and $\varepsilon = 0.01$ after $10^5$ iterations. 
\begin{figure}[ht]
\label{figure2}
\begin{center}
\rotatebox{0}{\resizebox{!}{0.3\linewidth}{%
   \includegraphics{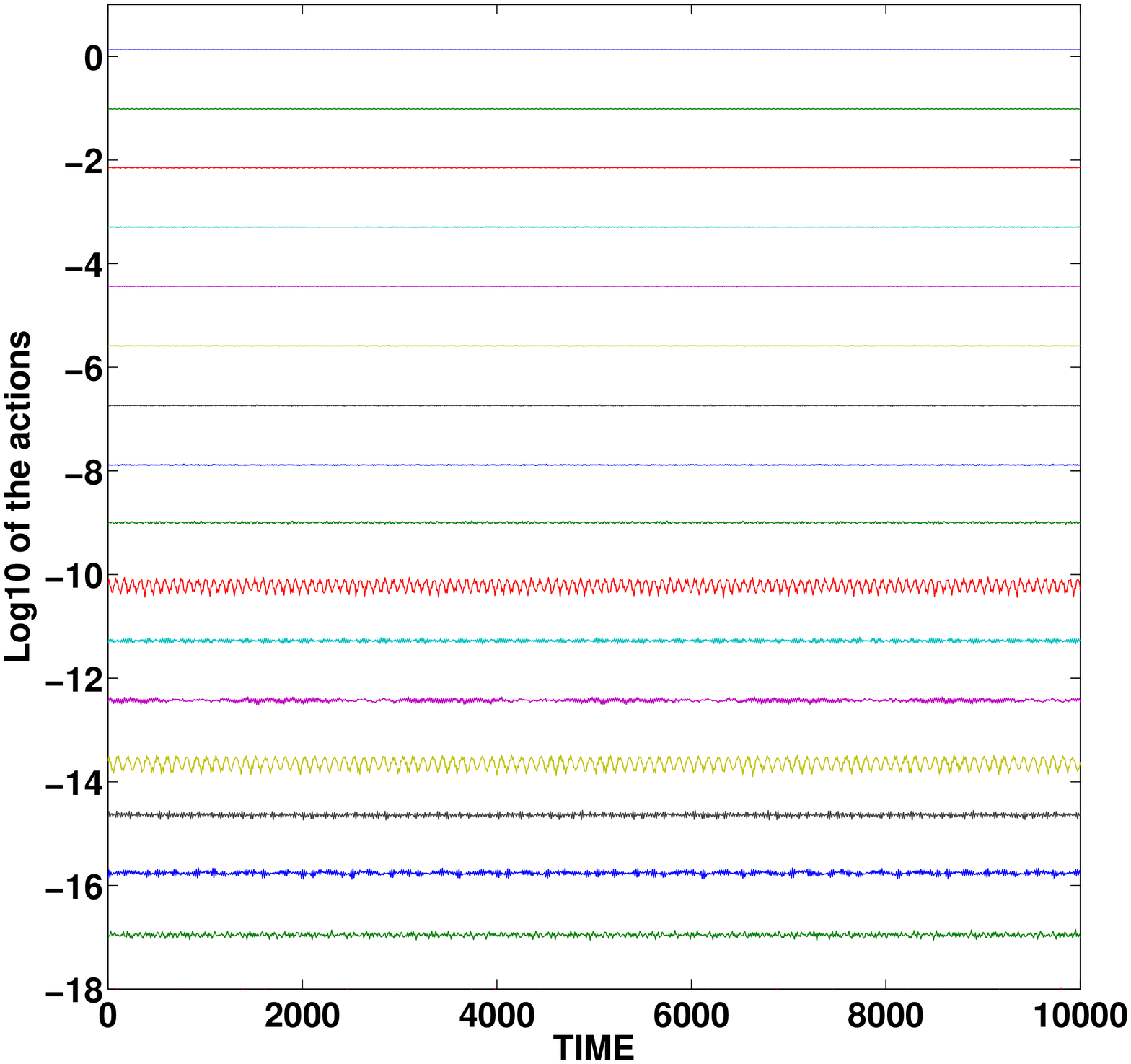}}} 
\rotatebox{0}{\resizebox{!}{0.3\linewidth}{%
   \includegraphics{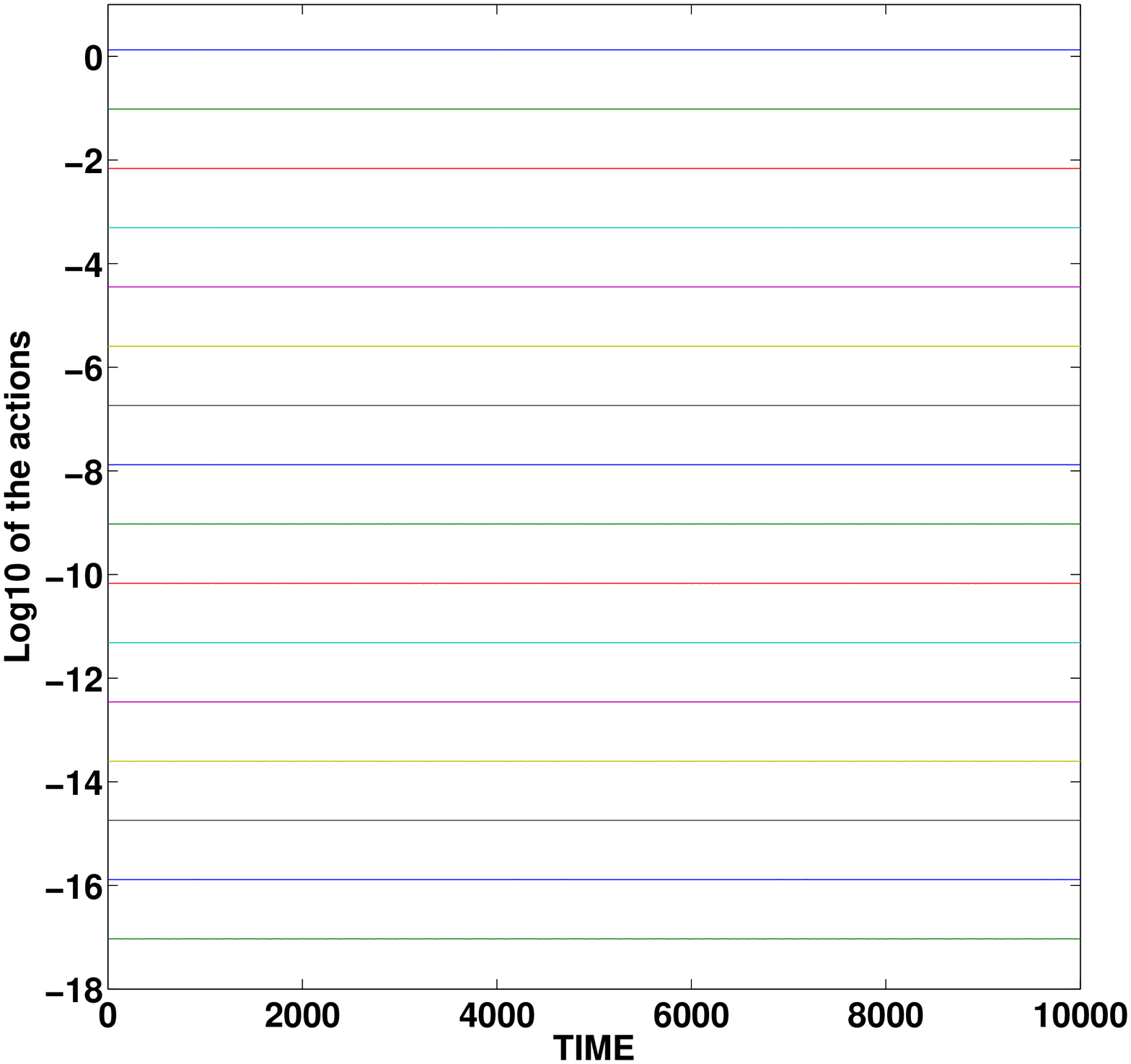}}} 
\caption{\small{Conservation of the actions for $\varepsilon = 0.1$ (left) and $\varepsilon = 0.01$ (right).} }
\end{center}
\end{figure}

\subsubsection{Wave equation on the circle}

We consider the wave equation on the circle
$$
u_{tt} - u_{xx} + m u = g(u), \quad x \in \T^1, \quad t \in \R, 
$$
where $m$ is a non negative real constant and $g$ a smooth real valued  function.
Introducing the variable $ v = u_t$, the corresponding Hamiltonian can be written
$$
H(u,v) =\int_{\T} \frac12 (v^2 +  u_x^2 +  m u^2) + G(u) \, \dd x,
$$
where $G$ is such that $\partial_u G = g$. Let $A := (-\partial_{xx} + m)^{1/2}$, and define the variables $(p,q)$ by
$$
q := A^{1/2} u, \quad\mbox{and}\quad p = A^{-1/2}v. 
$$
Then the Hamiltonian can be written 
$$
H = \frac12 \big(\langle Ap,p\rangle_{L^2} + \langle Aq,q\rangle_{L^2} \big) + \int_{\T} G(A^{-{1/2}}q) \, \dd x. 
$$
Let $\omega_a=\sqrt{|a|^2+m}$, $a \in \N =: \Nc$ be the eigenvalues of the operator $A$, and $\phi_a$ the associated eigenfunctions. Plugging the decompositions 
$$
q(x) = \sum_{a \in \N} q_a \phi_a(x) \quad \mbox{and} \quad p(x) = \sum_{a\in \N} p_a \phi_a(x)
$$
into the Hamiltonian functional, we see that it takes the form 
$$
H = \sum_{a \in \N}\omega_a \frac{p_a^2 + q_a^2}{2} + P
$$
where $P$ is a function of the variables $p_a$ and $q_a$. 
Using the complex coordinates 
$$
\xi_a = \frac{1}{\sqrt{2}} ( q_a + ip_a) \quad\mbox{and}\quad 
\eta_a = \frac{1}{\sqrt{2}} ( q_a - ip_a)
$$
the Hamiltonian function can be written under the form \eqref{E2} with a nonlinearity depending on $G$. As in the previous case, it can be shown  that the condition \eqref{nonres} is fulfilled for a set of constant $m$ of full measure (see \cite{BG06, Bam07}).   A collocation discretization on equidistant points of $[0,2\pi]$ yields the same discretization as previously (with $d = 1$). 

In this situation, the symmetric Strang splitting scheme 
$$
\varphi_{P^{(K)}}^{h/2} \circ \varphi_{H_0^{(K)}}^h \circ \varphi_{P^{(K)}}^{h/2}
$$
corresponds to the Deuflhard's method \cite{Deuf79}. If moreover we consider the Hamiltonian 
$$
H^{(K)}(z) = H_0^{(K)}(z) + P^{(K)}(\Phi(h \Omega) z)
$$
where $\Omega$ is the matrix with elements $\omega_a$, $a \in \Nc_K$, and $\Phi(x)$  a smooth function that is real, bounded and such that $\Phi(0) = 1$, then 
the splitting schemes associated with this decomposition coincide with the symplectic mollified impulse methods (see \cite[Chap. XIII]{HLW} and \cite{CHL08c}). 

\section{Proof of the normal form result}

\label{SProof}

The rest of the paper consists in proving Theorem \ref{TNF}. 

In the following, we denote by $\Tc_r$ the set of polynomial of order $r$ on $\C^{\Zc_K}$ (for sake of simplicity, we do note write the dependance in $K$ in the notation $\Tc_r$). If 
$$
Q = \sum_{\ell = 0}^{r} \sum_{\jb \in \Zc_K^\ell} Q_\jb z_{\jb}
$$
is an element of $\Tc_r$, we set
$$
\SNorm{Q}{\Tc_r} = \max_{\ell = 0,\ldots,r} \max_{\jb \in \Zc_K^\ell} | Q_\jb |.
$$
If moreover $Q \in \mathcal{C}([0,1],\Tc_r)$ we set
$$
\Norm{Q}{\Tc_r} = \max_{\lambda \in [0,1]} \SNorm{Q(\lambda)}{\Tc_r}. 
$$

Using the assumptions on $P^{(K)}$, we can write a Taylor expansion of $P$ around $0$, 
$$
P^{(K)}(z) =  P_r + Q_r = \sum_{\ell = 3}^r \sum_{\jb \in \Zc_K^\ell} P_\jb z_\jb + Q_r(z)
$$
where 
$$
|P_\jb| \leq C K^{\beta_0}
$$
where $C$ and $\beta_0$ depend on $\beta(\ell)$ and $C(\ell)$, $\ell = 0,\ldots,r$ in \eqref{EestP}. 

Notice that $Q_r(z) \in \mathcal{C}^{\infty}(\C^{\Zc_K},\C)$ admits of zero of order $r +1$ and satisfies 
$$
\Norm{X_{Q_r}(z)}{} \leq C K^{\beta_0} \Norm{z}{}^r
$$
for $z \in \Uc_K$, provided $\beta_0 = \beta_0(r,d)$ is large enough.

Before giving the proof of Theorem \ref{Tmain}, we give easy results on the flow of non autonomous polynomials Hamiltonian.
\begin{lemma}
\label{LP1}
Let $k \geq 1$ and let $P(\lambda) \in \mathcal{C}([0,1],\Tc_{k+1})$ be a homogeneous polynomial of order $k+1$ depending on $\lambda \in [0,1]$. Then
\begin{itemize}
\item[(i)]
There exists a constant $C$ depending on $k$ such that 
for all $z \in \C^{\Zc_K}$ and all $\lambda \in [0,1]$, we have
$$
|P(\lambda, z)| \leq CK^{d(k+1)}\Norm{P}{\Tc_{k+1}}\Norm{z}{}^{k+1}. 
$$
\item[(ii)]
There exists a constant $C$ depending on $k$ such that
for any $z \in \C^{\Zc_K}$ and all $\lambda \in [0,1]$, 
$$
\Norm{X_{P(\lambda)}(z)}{} \leq CK^{d(k+1)} \Norm{P}{\Tc_{k+1}} \Norm{z}{}^{k}. 
$$
\end{itemize}
Moreover, 
Let $k_1$ and $k_2$ two fixed integers. Let $P$ and $Q$ two homogeneous polynomials of degree $k_1+1$ and $k_2+1$ such that $P \in \mathcal{C}([0,1],\Tc_{k_1 + 1})$ and $Q \in \mathcal{C}([0,1],\Tc_{k_2 + 1})$. 
Then $\{P,Q\} \in \mathcal{C}([0,1],\Tc_{k_1 + k_2})$ and we have
$$
\Norm{\{P,Q\}}{\Tc_{k_1 + k_1}} \leq C \Norm{P}{\Tc_{k_1 + 1}} \Norm{Q}{\Tc_{k_2 + 1}}
$$
for some constant $C$ depending on $k_1$ and $k_2$. 
\end{lemma}
\begin{Proof}
We have 
$$
|P(\lambda, z)| \leq \Norm{P}{\Tc_{k+1}} \sum_{\jb \in \Zc_K^{k+1}} |z_\jb|
$$
where we have set  for $\jb = (j_1,\ldots,j_\ell) \in \Zc_K^\ell$, 
$$
|z_\jb| = |z_{j_1}| \cdots | z_{j_\ell}|. 
$$
Using $|z_j| \leq \Norm{z}{}$
we easily obtain {\em (i)} using $\sharp \Zc_K \leq (2K +1)^{d}$. The second statement is proven similarly. The estimate on the Poisson brackets is trivial. 
\end{Proof}

\begin{lemma}
\label{LP2}
Let $r \geq 3$, 
$$
Q(\lambda,z) = \sum_{\ell = 3}^r \sum_{\jb \in \Zc_K^\ell} Q_{\jb}(\lambda) z_{\jb}
$$
be an element of $\mathcal{C}([0,1],\Tc_{r})$. Let $\varphi_{Q(\lambda)}^\lambda$ be the flow associated with the non autonomous real Hamiltonian $Q(\lambda)$. Then there exist a constant $C_r$ depending on $r$ such that 
\begin{equation}
\label{Eball}
\rho < \mathrm{inf}\big(1/2, C_r  K^{-dr}\Norm{Q}{\Tc_{r}}^{-1}\big)\quad\Longrightarrow\quad
\forall\, \lambda \in [0,1], \quad \varphi_{Q(\lambda)}^\lambda(B_K(\rho)) \subset B_K(2\rho).
\end{equation}
Moreover, if $F(\lambda) \in \mathcal{C}([0,1], \mathcal{C}^\infty(B_K(\rho),\C))$ has a zero of order $r$ at the origin, then 
$F(\lambda) \circ \varphi_{Q(\lambda)}^\lambda$ has a zero of order $r$ at the origin in $B_K(\rho)$. 
\end{lemma}
\begin{Proof}
Let $z^\lambda =\varphi_{Q(\lambda)}^\lambda(z^0)$. Using the estimates of the previous lemma, we have 
$$
\begin{array}{rcl}
\displaystyle
\frac{\dd}{\dd \lambda} \Norm{z^\lambda}{}^2 &=& 2 \langle z^\lambda , X_{Q(\lambda)}(z^\lambda) \rangle \\[1ex]
&\leq&  c_r K^{dr}\Norm{Q}{\Tc_{r}}\Norm{z^\lambda}{}\Big(\Norm{z^\lambda}{}^2 + \Norm{z^\lambda}{}^{r-1}\Big)
\end{array}
$$
for some constant $c_r$ depending on $r$. Hence, as long as $\Norm{z^\lambda}{} \leq 1$, we have
$$
\displaystyle\frac{\dd}{\dd \lambda} \Norm{z^\lambda}{}^2 \leq 2 c_r K^{dr}\Norm{Q}{\Tc_{r}}\Norm{z^\lambda}{}^3. 
$$
By a standard comparison argument, we easily get that for $z^0 \in B_K(\rho)$ we have 
$$
\forall\, \lambda \in [0,1],\quad \Norm{z^\lambda}{} \leq 2 \Norm{z^0}{}. 
$$
This shows \eqref{Eball} and the rest follows. 
\end{Proof}

We give now the general strategy of the proof of the normal form Theorem \ref{TNF},  showing in particular the need of working with non autonomous Hamiltonians and of considering the non resonance condition \eqref{nonres2}. 

We consider  a fixed step size $h$ satisfying \eqref{nonres2}. As in this section $K$ will be considered as fixed, we denote shortly $P^{(K)}$ by $P$ and $H_0^{(K)}$ by $H_0$. 
We consider the propagator
$$
\varphi_{H_0}^h \circ \varphi_{P}^h = \varphi_{H_0}^h \circ \varphi_{hP}^1. 
$$
We embed this application into the family of applications
$$
\varphi_{H_0}^h \circ \varphi_{hP}^\lambda, \quad \lambda \in [0,1]. 
$$
Formally, we would like to find a real Hamiltonian $\chi = \chi(\lambda)$ and a real Hamiltonian under  normal form $Z = Z(\lambda)$ and such that 
\begin{equation}
\label{eq:flots}
\forall\, \lambda \in [0,1]\quad \varphi_{H_0}^h \circ \varphi_{hP}^\lambda \circ \varphi_{\chi(\lambda)}^\lambda =  \varphi_{\chi(\lambda)}^\lambda \circ\varphi_{H_0}^h \circ \varphi_{hZ(\lambda)}^\lambda. 
\end{equation}
Let $z^0 \in \C^{\Zc_K}$ and $z^\lambda = \varphi_{H_0}^h \circ \varphi_{hP}^\lambda \circ \varphi_{\chi(\lambda)}^\lambda (z^0)$. Deriving the previous equation with respect to $\lambda$ yields
\begin{multline*}
\frac{\dd z^\lambda}{ \dd \lambda} = 
(D_z \varphi_{H_0}^h)_{\varphi_{H_0}^{-h}(z^\lambda)} X_{hP}( \varphi_{H_0}^{-h}(z^\lambda) ) + \\[1ex]
(D_z( \varphi_{H_0}^h \circ \varphi_{hP}^\lambda ))_{ \varphi_{hP}^{-\lambda} \circ \varphi_{H_0}^{-h}(z^\lambda)}
X_{\chi(\lambda) } (\varphi_{hP}^{-\lambda} \circ \varphi_{H_0}^{-h}(z^\lambda)) . 
\end{multline*}
Using Lemma \ref{Lchange} that remains obviously valid for non autonomous Hamiltonian, we thus have 
$$
\frac{\dd z^\lambda}{ \dd \lambda} = X_{A(\lambda)}(z^\lambda)
$$
where $A(\lambda)$ it the time dependent real Hamiltonian given by 
$$
A(\lambda) = hP \circ \varphi_{H_0}^{-h}+ \chi(\lambda) \circ \varphi_{hP}^{-\lambda} \circ \varphi_{H_0}^{-h}. 
$$
Using the same calculations for the right-hand side, \eqref{eq:flots} is formally equivalent to the following equation (up to an integration constant)
\begin{equation}
\label{eq:tg}
\forall\, \lambda \in [0,1]\quad 
hP \circ \varphi_{H_0}^{-h} + \chi(\lambda) \circ \varphi_{hP}^{-\lambda} \circ \varphi_{H_0}^{-h} = \chi(\lambda) + hZ(\lambda) \circ \varphi^{-\lambda}_{\chi(\lambda)} \circ \varphi_{H_0}^{-h}. 
\end{equation}
which is equivalent to 
\begin{equation}
\label{eq:tg1}
\forall\, \lambda \in [0,1]\quad 
  \chi(\lambda) \circ \varphi_{H_0}^h - \chi(\lambda) \circ \varphi_{hP}^{-\lambda}= hP - hZ(\lambda) \circ \varphi^{-\lambda}_{\chi(\lambda)}   . 
\end{equation}
In the following, we will solve this equation in $\chi(\lambda)$ and $Z(\lambda)$ with a remainder term of order $r+1$ in $z$. 
So instead of \eqref{eq:tg1}, we will solve  the equation
\begin{equation}
\label{eq:tg2}
\forall\, \lambda \in [0,1]\quad 
  \chi(\lambda) \circ \varphi_{H_0}^h - \chi(\lambda) \circ \varphi_{hP}^{-\lambda}= hP - (h Z(\lambda) + R(\lambda)) \circ \varphi^{-\lambda}_{\chi(\lambda)} . 
\end{equation}
where the unknown are $\chi(\lambda)$, and $Z(\lambda)$ are polynomials of order $r$, with $Z$ under normal form, and where $R(\lambda)$ possesses a zero of order $r+1$ at the origin.

In the following, we formally write 
$$
\chi(\lambda) = \sum_{\ell = 3}^r \chi_{[\ell]}(\lambda)  :=  \sum_{\ell = 3}^r\sum_{\jb \in \Zc_K^\ell} \chi_\jb(\lambda)  z_{\jb}
$$
and 
$$
Z(\lambda) = \sum_{\ell = 3}^r Z_{[\ell]}(\lambda)  :=  \sum_{\ell = 3}^r\sum_{\jb \in \Zc_K^\ell} Z_\jb(\lambda)  z_{\jb}
$$
where here the coefficients $Z_\jb(\lambda)$ are unknown and where we denote by $\chi_{[\ell]}(\lambda)$ and $Z_{[\ell]}(\lambda)$ the homogeneous part of degree $\ell$ in the polynomials $\chi(\lambda)$ and $Z(\lambda)$. 

Identifying the coefficients of degree $\ell \leq r$ in  equation \eqref{eq:tg2}, we obtain
$$
\chi_{[\ell]}(\lambda) \circ \varphi_{H_0}^h - \chi_{[\ell]}(\lambda) = hP_{[\ell]} - h Z_{[\ell]}(\lambda) + h G_{[\ell]}(\lambda;\chi_*,P_*,Z_*) . 
$$
where $G$ is a real Hamiltonian homogeneous of degree $\ell$ depending on the polynomials $\chi_{[k]}$, $P_{[k]}$ and $Z_{[k]}$ for $k < \ell$. In particular, its coefficients are polynomial of order $\leq \ell$ of the coefficients $\chi_{j}$,  $P_{j}$ and $Z_{j}$ for $j \in \Zc_K^k$, $k < \ell$.  

Writing down the coefficients, this equation is equivalent to 
$$
\forall\, \jb \in \Zc_K^r \quad
( e^{ih \Omega(\jb)} - 1)\chi_{\jb} = hP_{\jb} - h Z_{\jb} + h G_{\jb}
$$
and hence we see that the key is to control the small divisors $ e^{ih \Omega(\jb)} - 1$ to solve these equations recursively. 


\begin{lemma}
Let $\chi(\lambda)$ 
be an element of $\mathcal{C}([0,1],\Tc_{r})$. Let $\tau(\lambda) := \varphi_{\chi(\lambda)}^\lambda$ be the flow associated with the non autonomous real Hamiltonian $\chi(\lambda)$. 
Let $g \in \mathcal{C}([0,1],\Tc_r)$, then we can write for all $\sigma_0 \in [0,1]$, 
\begin{multline}
\label{Ereccomm}
g(\sigma_0) \circ \tau(\sigma_0) = g(\sigma_0) \\[1ex]
+\sum_{k = 0}^{r-1} \int_{0}^{\sigma_0} \cdots \int_{0}^{\sigma_k} \Big( \mathrm{Ad}_{\chi(\sigma_k)}\circ \cdots \circ  \mathrm{Ad}_{\chi(\sigma_1)} g(\sigma_0)\Big) \dd \sigma_1 \cdots \dd \sigma_k  + R(\sigma_0) 
\end{multline}
where by definition $\mathrm{Ad}_P(Q) = \{Q,P\}$
\begin{equation}
\label{ERlambda}
R(\sigma_0) =  \int_{0}^{\sigma_0} \cdots \int_{0}^{\sigma_{r}} \Big( \mathrm{Ad}_{\chi(\sigma_{r})}\circ \cdots \circ  \mathrm{Ad}_{\chi(\sigma_1)} g(\sigma_0)\Big)\circ { \tau(\sigma_{r})}\,  \dd \sigma_1 \cdots \dd \sigma_{r}. 
\end{equation}
Each term in the sum in Eqn. \eqref{Ereccomm} belongs (at least) to the space $\mathcal{C}([0,1], \Tc_{kr})$.
The term $R(\sigma_0)$ defines an element of $\mathcal{C}([0,1],\mathcal{C}^\infty(\C^{\Zc_K},\C))$ and has a zero of order at least $r+1$ at the origin. 
\end{lemma}

The proof of this lemma is given in \cite{FGP2}. 

\medskip
As mentioned previously, 
for a given polynomial $\chi\in  \mathcal{C}([0,1], \Tc_r)$ with $r \geq 3$, we use the following notation 
\begin{equation}
\label{Echiell}
\chi(\lambda,z) = \sum_{\ell = 3}^r \chi_{[\ell]}(\lambda) = \sum_{\ell = 3}^r \sum_{\jb \in \Zc_K^\ell} \chi_{\jb}(\lambda) z_{\jb}
\end{equation}
where $\chi_{[\ell]}(\lambda) \in \mathcal{C}([0,1], \Tc_r)$ is a homogeneous polynomial of degree $\ell$.

\begin{proposition}
\label{Pcomposition}
Let $\chi(\lambda)$ 
be an element of $\mathcal{C}([0,1],\Tc_{r})$. Let $\varphi_{\chi(\lambda)}^\lambda$ be the flow associated with the non autonomous real Hamiltonian $\chi(\lambda)$. 
Let $g \in \mathcal{C}([0,1],\Tc_r)$, then we can write for all $\lambda \in [0,1]$, 
$$
g(\lambda) \circ \varphi^\lambda_{\chi(\lambda)} = S^{(r)}(\lambda) + T^{(r)}(\lambda)
$$
where 
\begin{itemize}
\item $S^{(r)}(\lambda) \in \mathcal{C}([0,1],\Tc_r)$. Moreover, if we write 
$$
S(z) = \sum_{\ell = 3}^r S_{[\ell]}(\lambda)
$$
where $S_{[\ell]}(\lambda)$ is a homogeneous polynomial of degree $\ell$, then we have  for all $\ell = 3,\ldots,r$, 
$$
S_{[\ell]}(\lambda) = g_{[\ell]}(\lambda) + G_{[\ell]}(\lambda;\chi_*,g_*)
$$
where $G_{[\ell]}(\lambda;\chi_*,g_*)$ is a homogeneous polynomial depending on $\lambda$ and 
the coefficients $S_j$ are polynomials of order $ < \ell$ of the coefficients appearing in the decomposition of $g$ and $\chi$. Moreover, we have
\begin{equation}
\label{Ebound1}
\Norm{G_{[\ell]}(\lambda;\chi_*,g_*)}{} \leq \Big( 1 + \sum_{m = 3}^{\ell - 1}  \Norm{g_{[m]}}{}^\ell\Big)\Big(1 + \sum_{m = 3}^{\ell - 1} \Norm{\chi_{[m]}}{}^\ell\Big).
\end{equation}
\item $T^{(r)}(\lambda)\in \mathcal{C}([0,1],\mathcal{C}^{\infty}(\C^{\Zc_K}, \C))$ has a zero of order at least $r +1$ at the origin and satisfies for all $z \in B_K(1/2)$, 
$$
\Norm{X_{T^{(r)}(\lambda)}(z)}{}\leq  C_r K^{2rd}  C_r(\chi_*,g_*) \Norm{z}{}^r
$$
where 
\begin{equation}
\label{Ebound2}
C_r(\chi_*,g_*) \leq C \Big( 1 + \sum_{m = 3}^{r}  \Norm{g_{[m]}}{\Tc_r}^r\Big)\Big(1 + \sum_{m = 3}^{r} \Norm{\chi_{[m]}}{\Tc_r}^r\Big)
\end{equation}
with $C$ depending on $r$.  
\end{itemize} 
\end{proposition}
\begin{Proof}
Using the previous lemma, we define $S^{(r)}$ as the polynomial part of degree less or equal to $r$ in the expression \eqref{Ereccomm}: this polynomial part may be computed iteratively, from the homogeneity degree 3 to $r$. Actually, every Poisson bracket appearing in \eqref{Ereccomm} is taken with a polynomial $\chi(\sigma_k)$, which decomposes into homogeneous polynomials with degree 3 at least. The terms appearing in the sum in \eqref{Ereccomm} hence have an increasing valuation, and this allows the iterative computation. The remainder terms, together with the term $R(\lambda)$ in \eqref{ERlambda}, define the term $T^{(r)}$ (which is an element of $\mathcal{C}([0,1],\Tc_{2r})$). The properties of $S^{(r)}(\lambda)$ and $T^{(r)}(\lambda)$ are then easily shown using Lemma \ref{LP1}. 
 \end{Proof}

The next result (Proposition \ref{PNF} below) yields the construction of the {\em normal form} term $\psi_K$ of Theorem \ref{TNF}. 

\begin{definition}
A polynomial $Z$ on $\C^{\Zc_K}$ is said to be in normal form if we can write it
$$
Z = \sum_{\ell = 3}^{r} \sum_{\jb \in \Ac_K^\ell} Z_\jb z_{\jb}. 
$$
where $\Ac_K^\ell$ is defined in the beginning of Subsection \ref{SSak}. 
\end{definition}

\begin{proposition}
\label{PNF}
Assume that $H := H^{(K)}$ satisfies \eqref{Edecomp} with $P := P^{(K)}$ fulfilling  \eqref{EestP} and assume that $h \leq h_0$ satisfies the hypothesis \eqref{nonres2}.  
 Then 
there exist 
\begin{itemize}
\item a polynomial $\chi \in \mathcal{C}([0,1], \Tc_r)$ 
$$
\chi(\lambda) = \sum_{\ell = 3}^r \chi_{[\ell]}(\lambda)  :=  \sum_{\ell = 3}^r\sum_{\jb \in \Zc_K^\ell} \chi_\jb(\lambda)  z_{\jb}
$$

\item 
a polynomial $Z \in \mathcal{C}([0,1], \Tc_r)$
$$
Z(\lambda) = \sum_{\ell = 3}^r Z_{[\ell]}(\lambda) := \sum_{\ell = 3}^r \sum_{\jb \in \Ac_K^\ell} Z_j(\lambda) z_{\jb}
$$
in normal form,
\item 
a function $R(\lambda)  \in \mathcal{C}([0,1],\mathcal{C}^\infty(B_K(\rho),\C))$ with $\rho < c_0 K^{-\beta}$ for some constant $c_0>0$ and $\beta>1$ depending on $r$ and $d$, and having a zero of order at least $r+1$ at the origin 
\end{itemize}
 such that the following equation holds: 
\begin{equation}
\label{eq:tg22}
\forall\, \lambda \in [0,1]\quad 
  \chi(\lambda) \circ \varphi_{H_0}^h - \chi(\lambda) \circ \varphi_{hP}^{-\lambda}= hP - (h Z(\lambda) + R(\lambda)) \circ \varphi^{-\lambda}_{\chi(\lambda)} . 
\end{equation}
Furthermore there exists a constant $C_0$ depending on $r$ and $d$ such that 
$$
\Norm{\chi}{\Tc_r} + \Norm{Z}{\Tc_r} \leq C_0 K^{\beta}
$$
and such that for all $\rho < c_0 K^{-\beta}$ and all $z \in B_K(\rho)$, we have
$$
\forall\, \lambda \in [0,1],\quad \Norm{X_{R(\lambda)}(z)}{} \leq C_0 K^{\beta} \Norm{z}{}^r. 
$$
\end{proposition}
\begin{Proof}
Identifying the coefficients of degree $\ell \leq r$ in the equation \eqref{eq:tg22}, we get
$$
\chi_{[\ell]} \circ \varphi_{H_0}^h - \chi_{[\ell]} = hP_{[\ell]} - h Z_{[\ell]} + h G_{[\ell]}(\chi_*,P_*,Z_*) . 
$$
where $G$ is a real Hamiltonian homogeneous of degree $\ell$ depending on the polynomials $\chi_{[k]}$, $P_{[k]}$ and $Z_{[k]}$ for $k < \ell$. In particular, its coefficients are polynomial of order $\leq \ell$ of the coefficients $\chi_{\jb}$,  $P_{\jb}$ and $Z_{\jb}$ for $\jb \in \Zc_K^k$, $k < \ell$ and satisfy bounds like \eqref{Ebound1}. 
Writing down the coefficients, this equation is equivalent to 
$$
\forall\, \jb \in \Ic_r \quad
( e^{ih \Omega(\jb)} - 1)\chi_{\jb} = hP_{\jb} - h Z_{\jb} + h G_{\jb}. 
$$
We solve this equation by setting
$$
Z_\jb = P_\jb + G_\jb\quad \mbox{and}\quad \chi_\jb = 0 \quad \mbox{for}\quad \jb \in \Ac_K^\ell
$$
and
$$
Z_\jb = 0 \quad \mbox{and} \quad \chi_{\jb} = \frac{h}{e^{ih \Omega(\jb)} - 1} (P_\jb +  G_\jb) \quad \mbox{for}\quad \jb \notin \Ac_K^\ell. 
$$
Using \eqref{nonres2} and the result of Proposition \ref{Pcomposition} we  get the claimed bound for some $\beta$ depending on $r$. 

To define $R$, we simply define it by the equation \eqref{eq:tg2}. By construction and the assumption on $P = P^{(K)}$, and using bounds of the form \eqref{Ebound2}, it is easy to show that it  satisfies the hypothesis. 
\end{Proof}

\begin{Proofof}{Theorem \ref{TNF}}
Integrating the equation \eqref{eq:tg2} in $\lambda$, it is clear that the following equation holds: 
$$
\forall\, \lambda \in [0,1]\quad \varphi_{H_0}^h \circ \varphi_{hP}^\lambda \circ \varphi_{\chi(\lambda)}^\lambda =  \varphi_{\chi(\lambda)}^\lambda \circ\varphi_{H_0}^h \circ \varphi_{hZ(\lambda) + R(\lambda)}^\lambda.
$$
Note that using Proposition \ref{Pcomposition} and \eqref{Eball} we  show that  for $z \in B_K(\rho)$ with  $\rho = cK^{-\beta}$  we have
$$
\Norm{\varphi^\lambda_{\chi(\lambda)}(z) - z}{} \leq C K^{\beta} \Norm{z}{}^2.
$$
This implies in particular that 
$$
\Norm{z}{} \leq  \Norm{\varphi^\lambda_{\chi(\lambda)}(z)}{} + CK^{-\beta} \Norm{z}{}
$$
For $K$ sufficiently large, this shows that $\varphi_{\chi(\lambda)}^\lambda$ is invertible and send $B_{K}(\rho)$ to $B_K(2\rho)$. 
Moreover, we have the estimate, for all $\lambda \in [0,1]$, 
$$
\Norm{\big(\varphi^\lambda_{\chi(\lambda)}\big)^{-1}(z) - z}{} \leq C K^{\beta} \Norm{z}{}^2.
$$
We then define $\tau_K = \varphi^1_{\chi(\lambda)}$ and $\psi_K = \varphi^1_{hZ(\lambda) + R(\lambda)}$ and verify that these applications satisfy the condition of the theorem for suitable constant $C$ and $\beta$. 
\end{Proofof}

\end{document}